\documentclass {article}
\usepackage{authblk}
\usepackage{etex}
\usepackage[utf8]{inputenc}
\usepackage[english]{babel}
\usepackage{amsmath}
\usepackage{amsthm}
\usepackage{amssymb}
\usepackage{sectsty}
\usepackage{titlesec}
\usepackage{color}
\usepackage[color,matrix,arrow]{xy}
\usepackage{amsgen}
\usepackage{amstext}
\usepackage{amsbsy}
\usepackage{amsopn}
\usepackage{amsfonts}
\usepackage{eepic}
\usepackage{graphicx}
\usepackage{epsf}
\usepackage{faktor}
\usepackage{float}
\usepackage{pstricks}
\usepackage{enumerate}
\usepackage{eqnarray}
\usepackage{faktor}
\xyoption{all}

\usepackage{pgf}
\usepackage{tikz-cd}
\usetikzlibrary{automata, arrows.meta, positioning,decorations.pathmorphing}

\newcommand{\footrecall}[1]{%
} 
\usetikzlibrary{snakes,shapes,arrows,automata}

\titleformat*{\section}{\large\bfseries}
\titleformat*{\subsection}{\normalsize \bfseries}

\newcommand{\N}{\mathbb{N}}
\newcommand{\Z}{\mathbb{Z}}

\newcommand{\End}{\text{End}}
\newcommand{\Ker}{\text{Ker}}
\newcommand{\TCC}{\text{TCC}}
\newcommand{\BCC}{\text{BCC}}
\newcommand{\GL}{\text{GL}}
\newcommand{\au}{ \overrightarrow{u}}
\newcommand{\ax}{ \overrightarrow{x}}
\newcommand{\ay}{ \overrightarrow{y}}

\newcommand{\Inn}{\text{Inn}}

\newcommand{\GCP}{\text{GCP}}
\newcommand{\GTCP}{\text{GTCP}}
\newcommand{\GBrCP}{\text{GBrCP}}
\newcommand{\GBrP}{\text{GBrP}}

\newcommand{\Aut}{\text{Aut}}
\newcommand{\Via}{\text{Via}}
\newcommand{\Img}{\text{Im}}

\newcommand{\Rat}{\text{Rat}}

\newcommand{\Fin}{\text{Fin}}
\newcommand{\Coset}{\text{Coset}}
\newcommand{\Rec}{\text{Rec}}
\newcommand{\Alg}{\text{Alg}}
\newcommand{\CF}{\text{CF}}

\newcommand{\dphi}{\Delta_\phi}
\newcommand{\phitcc}{\phi\text{-}\TCC}
\newcommand{\psitcc}{\Psi\text{-}\TCC}
\newcommand{\Psitcc}{\psi\text{-}\TCC}
\newcommand{\phibcc}{\phi\text{-}\BCC}

\newcommand{\mc}{\mathcal}

\theoremstyle{definition}
\newtheorem{theorem}{Theorem}[section]
\newtheorem{corollary}[theorem]{Corollary}

\newtheorem{proposition}[theorem]{Proposition}
\newtheorem{question}[theorem]{Question}

\newtheorem{lemma}[theorem]{Lemma}

\newtheorem{remark}[theorem]{Remark}

\begin{document}

\title{On generalized (twisted) conjugacy separability of some extensions of groups}
\author{Andr\'e Carvalho\thanks{andrecruzcarvalho@gmail.com}}
\affil{Centre of Mathematics of the University of Porto

 R. Campo Alegre, 4169-007 Porto,
 
Portugal}
\maketitle

\begin{abstract}
We introduce separability properties corresponding to generalized versions of the conjugacy, twisted conjugacy, Brinkmann and Brinkmann's conjugacy problems and how they relate when finite and cyclic extensions of groups are taken. In particular, we prove that some (concrete) generalizations of  twisted conjugacy separability of a group $G$ with respect to virtually inner automorphisms are equivalent to some (concrete) generalizations of the conjugacy problem in finite extensions of $G$. Similarly, (generalized) conjugacy separability in cyclic extensions of $G$ implies (generalized) twisted conjugacy and Brinkmann's conjugacy separability in $G$. Applications include results in free, virtually abelian, virtually polycyclic groups and a proof that virtually free times free groups are conjugacy separable.
\end{abstract}
\section{Introduction}
 All groups  in this paper are assumed to be finitely presented.
The profinite topology on a group $G$ has as basis sets the cosets (left and right) of finite index subgroups of $G$. A group endowed with the profinite topology is a topological group.

Let $S$ be a subset of $G$. The following conditions are equivalent:
\begin{itemize}
\item $S$ is closed in the profinite topology
\item for each $g\not\in S$, there is some finite index normal subgroup $N\trianglelefteq_{f.i.} G$ such that $g\not\in SN$
\item for each $g\not\in S$, there is some surjective homomorphism $\phi$ from $G$ to a finite group such that $g\phi\not\in S\phi$
\end{itemize}

So, a group is residually finite if singletons are separable. A subset satisfying any (thus all) of the above conditions is said to be \emph{separable}. We will often just call such a set \emph{closed}.

Separability can be important in studying algorithmic problems. Indeed, if a subset $S$ of a finitely presented group is recursively enumerable, separable and its image under a surjective homomorphism to a finite group can be computed, then we can decide membership on $S$: on input such a set $S$ and a word $w$ on the generators, we run two algorithms in parallel, one that will recursively enumerate elements in $S$ and test if they are equal to $w$ using the positive part of the word problem and one that will enumerate finite quotients of $G$, compute the image of $S$ in this quotient and check if the element represented by $w$ belongs (in the quotient) to $S$. If $w$ represents an element of $S$, the first algorithm will stop and we answer \texttt{YES}, while the second will not stop; if $w$ represents an element of $G\setminus S$, the second algorithm will stop and we answer \texttt{NO}, while the first will not stop. This argument was used by Mal'cev \cite{[Mal58],[Mal83]} to prove that residually finite groups have decidable word problem. Similar notions exist for other algorithmic problems: Mostowski used a similar argument to show that conjugacy separable groups, that is, groups having separable conjugacy classes, have decidable conjugacy problem in \cite{[Mos66]}; subgroup separable groups, also known as LERF - locally extended residually finite - have decidable subgroup membership problem; and twisted conjugacy separable groups have decidable twisted conjugacy problem \cite{[FT06],[FT07],[Ter24]}. Extensive work has been done in establishing separability properties for different classes of groups: for example, Ribes and Zalesskii \cite{[RZ93]} proved that free groups are product separable (meaning that a product $H_1\cdots H_n$ of finitely generated subgroups is separable), which, by the work of Pin and Reutenauer \cite{[PR91]}, means that the closure of a rational subset of a free group can be computed (and closed rational subsets are exactly the sets that can be written as finite unions of sets of the form $aH_1\cdots H_n$) and implies the Rhodes Type II Conjecture (see \cite{[Pin89]} for a survey on the problem); 
conjugacy separability has been proved for virtually free groups \cite{[Dye79]}, virtually polycyclic groups \cite{[For76],[Rem69]}, RAAGs \cite{[Min12]}, one-relator groups with torsion \cite{[MZ13]}, and others; twisted conjugacy separability is proved for polycyclic-by-finite groups in \cite{[FT06]} - further results were obtained in \cite{[Ter24]}; effective separability has been studied by many (see \cite{[DFP22]} for a survey); Gitik and Rips proved in \cite{[GR95]} that free groups are double coset separable, Mart\'inez-Pedroza and Sisto proved that double cosets of relatively quasiconvex subgroups with compatible parabolics in QCERF relatively hyperbolic groups were separable in \cite{[MS12]}.

In terms of decidability, it is well known that proving decidability in certain extensions of groups reduces  to proving  decidability of some other (\emph{harder}) algorithmic problem on the starting group: for example, in \cite{[BBV05b],[LS11],[Lev23]}, the twisted conjugacy problem with respect to virtually inner automorphisms (those having a power which is inner) of free groups is used to study the conjugacy problem in virtually free groups; in \cite{[BMMV06]}, Bogopolski, Martino, Maslakova and Ventura connect the conjugacy problem in free-by-cyclic groups with  Brinkmann's conjugacy problem and the twisted conjugacy problem in free groups; in \cite{[Log22]}, Logan shows that proving the decidability of the twisted conjugacy problem and of a two-sided version of Brinkmann's problem can be used to prove the decidability of the conjugacy problem in ascending HNN-extensions of groups and this was later used in \cite{[CD24c]} to prove the decidability of the conjugacy problem in ascending HNN-extensions of free-abelian times free groups.

Our goal is to establish separability analogues of the results concerning finite and cyclic extensions relating generalized versions of these problems (concrete definitions will appear in the next section). To do so, we start by introducing new (generalized) separability concepts and show that establishing them implies decidability of the corresponding (generalized) decision problems.

Then, we focus on  finite extensions. It is proved in \cite{[Gor86]} that conjugacy separability is not preserved by taking finite extensions and it is not known if twisted conjugacy separability is preserved by taking finite extensions (see \cite{[Ter24]}).  In \cite{[CZ09]}, an example of a non-conjugacy separable group with a finite index conjugacy separable subgroup is presented and a sufficient condition for this property to be preserved by taking finite index subgroups is shown.
We prove that depending on the type of generalization, conjugacy separability in a finite extension is equivalent to twisted conjugacy separability with respect to virtually inner automorphisms (the set of all such automorphisms of a group $G$ will be denoted by $\Via(G)$). In particular, in the simple case, we have the following result.

\newtheorem*{simple finite equiv}{Corollary \ref{simple finite equiv}}
\begin{simple finite equiv}
Let $\mc C$ be a class of finitely generated groups and $\mc{FEC}$ be the class of finite extensions of groups in $\mc C$. Then, groups in $\mc C$ are $\Via$-twisted conjugacy separable if and only if groups in $\mc{FEC}$ are conjugacy separable.
\end{simple finite equiv}

It is not known whether free groups are twisted conjugacy separable or not. Since finitely generated virtually free groups are conjugacy separable \cite{[Dye79],[Ste70]}, we get that that is the case if the automorphism is virtually inner:
\newtheorem*{free via twisted}{Corollary \ref{free via twisted}}
\begin{free via twisted}
Free groups are $\Via$-twisted conjugacy separable.
\end{free via twisted}

We can also used the result (in a different direction) to prove that virtually abelian groups are Coset-generalized conjugacy separable.
\newtheorem*{virtually abelian conjugacy}{Corollary \ref{virtually abelian conjugacy}}
\begin{virtually abelian conjugacy}
Every virtually abelian group is $\Coset$-generalized conjugacy separable.
\end{virtually abelian conjugacy}

Finally, we prove that groups containing a direct product of free groups as a finite index subgroup are conjugacy separable. This class includes virtually $F_n\times \Z$ groups, which contain  non-elementary unimodular generalized Baumslag-Solitar groups \cite{[Lev07]} and all the subgroup separable free-by-cyclic groups with a polynomially growing monodromy \cite{[Kud24]}.

\newtheorem*{main direct}{Theorem \ref{main direct}}
\begin{main direct}
Virtually free times free groups are conjugacy separable.
\end{main direct}

Regarding cyclic extensions, we can prove that if $G\rtimes \Z$ is (generalized) conjugacy separable, then $G$ is (generalized) twisted conjugacy separable and (generalized) Brinkmann's conjugacy separable.  For the simple versions of the problems, we have the following:

\newtheorem*{simple cyclic equiv}{Corollary \ref{simple cyclic equiv}}
\begin{simple cyclic equiv}
If the class of $G\rtimes \Z$ groups is conjugacy separable, then $G$ is twisted conjugacy separable and Brinkmann conjugacy separable.
\end{simple cyclic equiv}
The relation between conjugacy separability of $G\rtimes \Z$ and twisted conjugacy separability of $G$ has been proved in \cite{[FT06]}. Here, we not only prove the connection with Brinkmann conjugacy separability, but also prove this result for the generalized notions.

Since virtually polycyclic groups are conjugacy separable and cyclic extensions of virtually polycyclic groups are virtually polycyclic, it follows naturally that virtually polycyclic groups are twisted conjugacy separable (this was already proved in \cite{[FT06]}) and Brinkmann conjugacy separable.
 
 \newtheorem*{virt polycyclic twisted Brink}{Corollary \ref{virt polycyclic twisted Brink}}
\begin{virt polycyclic twisted Brink}
Virtually polycyclic groups are twisted conjugacy separable and Brinkmann conjugacy separable.
\end{virt polycyclic twisted Brink}

In Section 2, we introduce some preliminary notions on subsets of groups. In Section 3, we introduce the separability analogues of the generalized conjugacy problem, twisted conjugacy problem, Brinkmann's problem and Brinkmann's  conjugacy problem. In Section 4, we focus on finite extensions: we provide an alternative proof that finite extensions of product separable groups are product separable by describing the closure of rational subsets of finite extensions of groups using the closure of subsets of the original group and show that, for certain generalizations (which include the classical simple case), being twisted conjugacy separable with respect to virtually inner automorphisms is equivalent to conjugacy separability on finite extensions of the group. This implies in particular that free groups are twisted conjugacy separable with respect to virtually inner automorphisms and that virtually abelian groups are (Coset)-generalized conjugacy separable.  Finally, in Section 5, we turn our attention to cyclic extensions of groups and we prove that if $G\rtimes \Z$ groups are conjugacy separable, then $G$ is Brinkmann conjugacy separable and twisted conjugacy separable, which implies that virtually polycyclic groups are  twisted conjugacy separable and Brinkmann conjugacy separable.

\section{Preliminaries}
Here, we introduce some notions about subsets of groups defined by language-theoretic conditions.

Let $G=\langle A\rangle$ be a finitely generated group, $A$ be a finite generating set, $\widetilde A=A\cup A^{-1}$ and $\pi:\widetilde A^*\to G$ be the canonical (surjective) homomorphism. This notation will be kept throughout the paper.
A subset $K\subseteq G$ is said to be \emph{rational} if there is some rational language $L\subseteq \widetilde A^*$ such that $L\pi=K$ and \emph{recognizable} if $K\pi^{-1}$ is rational. 
We will denote by $\Rat(G)$ and $\Rec(G)$ the class of rational and recognizable subsets of $G$, respectively. Rational subsets generalize the notion of finitely generated subgroups.

\begin{theorem}[\cite{[Ber79]}, Theorem III.2.7]
\label{AnisimovSeifert}
Let $H$ be a subgroup of a group $G$. Then $H\in \Rat(G)$ if and only if $H$ is finitely generated.
\end{theorem}

Similarly, recognizable subsets generalize the notion of finite index subgroups.

\begin{proposition}
\label{rec fi}
Let $H$ be a subgroup of a group $G$. Then $H\in \Rec(G)$ if and only if $H$ has finite index in $G$.
\end{proposition}

In fact, if $G$ is a group and $K$ is a subset of $G$ then $K$ is recognizable if and only if $K$ is a (finite) union of cosets of a subgroup of finite index.

A natural generalization of these concepts concerns the class of context-free languages. 
A subset $K\subseteq G$ is said to be \emph{algebraic} if there is some context-free language $L\subseteq \widetilde A^*$ such that $L\pi=K$ and \emph{context-free} if $K\pi^{-1}$ is context-free. 
We will denote by $\Alg(G)$ and $\CF(G)$ the class of algebraic and context-free subsets of $G$, respectively.  It follows from \cite[Lemma 2.1]{[Her91]} that these definitions, as well as the definitions of rational and recognizable subsets, do not depend on the finite alphabet $A$ or the surjective homomorphism $\pi$.

It is obvious from the definitions that $\Rec(G)$, $\Rat(G)$, $\CF(G)$ and $\Alg(G)$ are closed under union, since both rational and context-free languages are closed under union. The intersection case is distinct: from the fact that rational languages are closed under intersection, it follows that $\Rec(G)$ must be closed under intersection too. However $\Rat(G)$, $\Alg(G)$ and $\CF(G)$ might not be. Another important closure property is given by the following lemma from \cite{[Her91]}.

\begin{lemma}\cite[Lemma 4.1]{[Her91]}
\label{herbst letra}
 Let $G$ be a finitely generated group, $R\in \Rat(G)$ and $C \in \{\Rec, \CF\}$.
If $K \in C(G)$, then $KR, RK \in C(G).$
\end{lemma}

This lemma is particularly useful in the case where $R$ is a singleton.

The following theorem was obtained by Grunschlag and Silva in the case of rational and recognizable subsets and by the author in the case of algebraic and context-free subsets.

\begin{theorem}[\cite{[Gru99],[Sil02b],[Car23c]}]\label{grunschlag-silva-carvalho}
Let $G$ be a finitely generated group, $K\in\{\Rat, \Rec, \Alg,\CF\}$ and $H\leq_{f.i.} G$. If $G$ is the disjoint union $G=\cup_{i=1}^n Hb_i$, then $K(G)$ consists of all subsets of the form 
\begin{align*}
\bigcup_{i=1}^{n} L_ib_i \qquad (L_i\in K(H)).
\end{align*}
\end{theorem}

Finally, we will denote the inner automorphism representing conjugation by an element $g\in G$ by $\lambda_g$, that is, for all $x\in G$, 
$x\lambda_g=g^{-1}xg$.

\section{Generalized separability}
In this section, we will continue the work mentioned above and introduce the  separability properties associated with the following algorithmic problems: let $\mc C$ be a class of subsets of a group $G$.

\begin{itemize}
\item \textbf{$\GCP_{\mc C}(G)$ -  $\mc C$-generalized conjugacy problem:} taking as input a subset $K\in \mc C$ and an element $x\in G$, decide whether there is a conjugate of $x$ in $K$;
\item \textbf{$\GTCP_{\mc C}(G)$ - $\mc C$-generalized twisted conjugacy problem:} taking as input a subset $K\in \mc C$, an automorphism $\phi\in \Aut(G)$ and  an element $x\in G$, decide whether $x$ has a $\phi$-twisted conjugate  in $K$, i.e., whether there is some $z\in G$ such that $(z^{-1}\phi) x z\in K$;
\item \textbf{$\GBrCP_{\mc C}(G)$  -  $\mc C$-generalized Brinkmann's conjugacy problem:} taking as input a subset $K\in \mc C$, an automorphism $\phi\in \Aut(G)$ and  an element $x\in G$, decide whether there is some $k\in \Z$ such that $x\phi^k$ has a conjugate in $K$;
\item \textbf{$\GBrP_{\mc C}(G)$ -  $\mc C$-generalized Brinkmann's problem:} taking as input a subset $K\in \mc C$, an automorphism $\phi\in \Aut(G)$ and an element $x\in G$, decide whether there is some $k\in \Z$ such that $x\phi^k$ belongs to $K$.
\end{itemize}

All these problems have been considered in \cite{[Car23b]} for decidability purposes. The generalized conjugacy problem was proven to be decidable for virtually free groups in \cite{[LS11]}. Brinkmann's conjugacy problem was proven to be decidable for automorphisms of the free group in \cite{[Bri10]} and this was later extended to general endomorphisms in \cite{[Log22],[CD24]}. Results on other classes of groups have also been obtained in \cite{[CD24b],[CD24c]}. This problem is important as it plays an important role in proving decidability of the conjugacy problem in free-by-cyclic groups (see \cite{[BMMV06]}).

Given a subset $S\subseteq G$, we will write $\alpha(S)$ to denote the union of all conjugacy classes of elements in $S$, that is $$\alpha(S)=\{x^{-1}sx\mid x\in G, s\in S\}.$$

For an automorphism $\phi\in \Aut(G)$, we will denote the union of all $\phi$-twisted conjugacy classes of elements in a set $S$ with a $\phi$-twisted conjugator in a set $N$ by $\phitcc_N(S)$, that is, 
$$\phitcc_N(S)=\{(x^{-1}\phi)sx\mid x\in N, s\in S\}.$$ When $N$ is the whole group, we will simply write $\phitcc(S)$. 

We denote the $\phi$-orbit of $x$ by $\mc O_\phi(x)$ and given a subset $S$ we write $\mc O_\phi(S)$ to denote $\bigcup_{s\in S}\mc O_\phi(s)$ . The union of \emph{Brinkmann conjugacy classes} will be denoted by $$\phibcc(S)=\{x^{-1}sx\mid x\in G, s\in \mc O_\phi(S)\}=\alpha(\mc O_\phi(K)).$$
Occasionally, we will also use a subscript on $\BCC$ and $\alpha$ to indicate which group the conjugators belong to.

We can now define a separability property for each of the algorithmic problems defined above. Given a class $\mc C$ of subsets, we say that a group is:

\begin{itemize}
\item \textbf{$\mc C$-generalized conjugacy separable} if $\alpha(K)$ is separable for all $K\in \mc C$;
\item \textbf{$\mc C$-generalized twisted conjugacy separable} if  $\phitcc(K)$ is separable for all  $K\in \mc C$ and $\phi\in \Aut(G)$;
\item \textbf{$\mc C$-generalized Brinkmann's conjugacy separable}  if  $\phibcc(K)$ is separable for all  $K\in \mc C$ and $\phi\in \Aut(G)$;
\item \textbf{$\mc C$-generalized Brinkmann's separable:} if $\mc O_\phi(K)$ is separable for all $K\in \mc C$ and $\phi\in \Aut(G)$.
\end{itemize}

The main classes of subsets that we will consider will be the classes of finite subsets, cosets of finitely generated subgroups, and the classes of rational, recognizable, algebraic and context-free subsets, and we will denote them by $\Fin$, $\Coset$, $\Rat$, $\Rec$, $\Alg$ and $\CF$, respectively. Notice that the $\Fin$ case is equivalent to the usual simple case, where sets are singletons because a finite set is a finite union of singletons and the finite union of closed sets is still closed.

If $N\leq_{f.i.} G$ is characteristic, then an automorphism $\phi\in \Aut(G)$ induces an automorphism $\overline\phi\in \Aut(\faktor{G}{N})$ such that $\pi\overline\phi=\phi\pi$, where $\pi:G\to \faktor{G}{N}$ is the standard projection. Notice that, then $\phi^k\pi=\pi\overline\phi^k$, for all $k\in \N$.
 \begin{figure}[H]
\begin{equation*}
\begin{tikzcd}[ampersand replacement=\&]
G \arrow{r}{\pi} \arrow[d,"\phi"']\& \faktor{G}{N} \arrow[d,dashed,"{\overline\phi}"] \\
G \arrow{r}{\pi} \& \faktor{G}{N}
\end{tikzcd}
\end{equation*}
 \end{figure}
We will keep this notation throughout the paper. We say that a subset $S\subseteq G$ is finite-computable if, given a homomorphism $\phi$ from $G$ to a finite group $F$, then $S\phi\subseteq F$ can be computed.

We remark that the subsets $\alpha(K)$, $\phitcc(K)$, $\phibcc(K)$ and $\mc O_\phi(K)$ are recursively enumerable if $K$ is recursively enumerable. It is easy to see that, if $\psi:G\to F$ is a surjective homomorphism to a finite group, then $(\alpha(K))\psi=\alpha(K\psi)$, which is computable, and so, if $\mc C$ is a class of finite-computable subsets,  $\mc C$-generalized conjugacy separability implies $\mc C$-generalized conjugacy decidability.

We will prove that the same holds for the remaining problems. But first, we will show a technical lemma. The proof is essentially the same as in  \cite[Proposition 4.0.2]{[Car23t]}, but we include it here for completeness, as the original statement refers only to virtually free groups. We say that a subgroup $H\leq G$ is fully invariant if $H\varphi\subseteq H$, for all $\varphi\in \End(G)$
\begin{lemma} \label{fully invariant}
Let $G$ be a group and $F\leq_{f.i.} G$ be such that membership in $N$ can be tested. Then, given $N$ and a transversal for $N$, we can compute a fully invariant finite index subgroup $N\leq_{f.i.} F$ and a right transversal for $N$.
\end{lemma}
\noindent\textit{Proof.}  By \cite[Lemma 4.0.1]{[Car23t]}, the intersection $N$ of all  normal subgroups of $G$ of index at most $[G:F]$  is a fully invariant finite index subgroup. Also, $N\leq F$ by construction. We will now prove that $N$ is computable. We start by  enumerating all finite groups of cardinality at most $[G:F]$. For each such group $K=\{k_1,\ldots, k_s\}$ we enumerate all homomorphisms from $G$ to $K$ by defining images of the generators and checking all the relations. For each homomorphism $\theta:G\to K$, we have that $[G:\Ker(\theta)]=|\Img(\theta)|\leq |K|\leq m$. In fact, all normal subgroups of $G$ of index at most $m$ are of this form. We compute generators for the kernel of each $\theta$, which is possible by Schreier's Lemma, since we can test membership in $\Ker(\theta)$, which is a finite index subgroup. We can also find $a_2,\ldots,a_{s}\in G$ such that 
$$G=\Ker(\theta)\cup\Ker(\theta)a_2\cup\cdots\cup \Ker(\theta)a_{s},$$
taking $a_i$ such that $a_i\theta=k_i.$

Hence, we can compute all finite subgroups of index at most $[G:F]$ and we can test membership in each of them (by checking if the image under the corresponding homomorphism $\theta$ is trivial). Thus, using Schreier's lemma we can compute the intersection $N$ of all such subgroups and 
, and a decomposition of $G$ as a disjoint union
\begin{align*}
G=Nb_1\cup Fb_2\cup \cdots \cup Nb_r.
\end{align*}
since membership in $N$ can be tested, as it amount to checking membership each finite index subgroup of index at most $[G:F]$.
\qed\\

\begin{theorem}
Let $G$ be a finitely presented group, $\phi\in \Aut(G)$,  and $\mc C$ be a class of  finite-computable subsets and let $\psi:G \to F$ be a surjective homomorphism to a finite group $F$. Then, for all $K\in \mc C$,  we can compute $\phitcc(K)\psi$, $\mc O_\phi(K)\psi$, and $\phibcc(K)\psi$.
\end{theorem}
\noindent\textit{Proof.} We have that $\Ker(\psi)\leq_{f.i.} G$ and, using Schreier's Lemma, we can compute a finite generating set and a transversal for $\Ker(\psi)$. 
By Lemma \ref{fully invariant}, we can compute generators and a transversal for a characteristic (in fact, for a fully invariant) subgroup $N\leq_{f.i.} \Ker(\psi)$ and compute the projection $\pi:G\to\faktor{G}{N}$.

We will prove that $\phitcc(K)\psi=\overline\phi
\text{-TCC}(K\pi)\pi^{-1}\psi$, $\mc O_\phi(K)\psi=\mc O_{\overline \phi}(K\pi)\pi^{-1}\psi$ and $\phibcc(K)\psi=\overline\phi
\text{-BCC}(K\pi)\pi^{-1}\psi$.

Given $k\in K$ and $x\in G$, we have that $$((x^{-1}\phi) k x)\pi=(x^{-1}\phi\pi) (k\pi) (x\pi)=(x^{-1}\pi\overline\phi)(k\pi) (x\pi),$$
so $\phitcc(K)\subseteq \overline\phi
\text{-TCC}(K\pi)\pi^{-1}$ and $\phitcc(K)\psi\subseteq \overline\phi
\text{-TCC}(K\pi)\pi^{-1}\psi$. 

Conversely, if $x\in \overline\phi
\text{-TCC}(K\pi)\pi^{-1}$, that is $x\pi\in  \overline\phi
\text{-TCC}(K\pi)$, then $x\pi= (z^{-1}\overline\phi)(k\pi)z$, for some $z\in \faktor{G}{N}$ and $k\in K$. Since $\pi$ is surjective, then $z=y\pi$ for some $y\in G$ and $$x\pi= (y^{-1}\pi\overline\phi)(k\pi)
(y\pi)=(y^{-1}\phi\pi)(k\pi)(y\pi)=((y^{-1}\phi)ky)\pi.$$
Since $N\leq \Ker(\psi)$, then $x\psi=((y^{-1}\phi)ky)\psi\in \phitcc(K)\psi$ and $\phitcc(K)\psi= \overline\phi
\text{-TCC}(K\pi)\pi^{-1}\psi$. 

We deal with the remaining cases similarly.
If $x\in \mc O_\phi(K)$, then $x=k\phi^n$ for some $k\in K$ and $n\in \Z$. Then $$x\pi=k\phi^n\pi=k\pi\overline\phi^n$$ and $x\in \mc O_{\overline \phi}(K\pi)\pi^{-1}$. Then,  $\mc O_\phi(K)\psi\subseteq \mc O_{\overline \phi}(K\pi)\pi^{-1}\psi$.

Conversely, if $x\in \mc O_{\overline \phi}(K\pi)\pi^{-1}$, that is, $x\pi\in  \mc O_{\overline \phi}(K\pi)$, then $x\pi=k\pi\overline\phi^n=k\phi^n\pi$ for some $k\in K$ and $n\in \Z$. Then, since $N\leq \Ker(\psi)$, $x\psi=k\phi^n\psi$ and $x\psi \in \mc O_\phi(K)\psi$. Hence,  $\mc O_\phi(K)\psi= \mc O_{\overline \phi}(K\pi)\pi^{-1}\psi$.

If $x\in \phibcc(K)$, then $x=z^{-1}(k\phi^n)z$ for some $z\in G$ and $n\in \Z$. Then 
$$x\pi=(z^{-1}(k\phi^n)z)\pi=(z^{-1}\pi)(k\phi^n\pi)(z\pi)=(z^{-1}\pi)(k\pi\overline\phi^n)(z\pi)\in \overline\phi
\text{-BCC}(K\pi).$$
Hence, $x\in \overline\phi
\text{-BCC}(K\pi)\pi^{-1}$, so $\phibcc(K)\psi\subseteq\overline\phi
\text{-BCC}(K\pi)\pi^{-1}\psi$.

Conversely, if $x\in \overline\phi\text{-BCC}(K\pi)\pi^{-1}$, that is, $x\pi\in \overline\phi\text{-BCC}(K\pi)$, then there are $k\in K$,  $n\in \Z$ and $z\in \faktor{G}{N}$ such that 
$x\pi=z^{-1}(k\pi\overline\phi^n)z.$ Since $\pi$ is surjective, $z=y\pi$ for some $y\in G$ and 
$$x\pi=(y^{-1}\pi)(k\pi\overline\phi^n)(y\pi)=(y^{-1}\pi)(k\phi^n\pi)(y\pi)=(y^{-1}(k\phi^n)y)\pi$$ and, since $N\leq \Ker(\psi)$, 
$x\psi=(y^{-1}(k\phi^n)y)\psi\in \phibcc(K)\psi$.

Therefore,  $\phibcc(K)\psi=
 \overline\phi\text{-BCC}(K\pi)\pi^{-1}\psi$.
 
 Now, it is clear that, since $K$ is finite-computable and $\faktor{G}{N}$ is finite and computable, then $\overline\phi
\text{-TCC}(K\pi)$, $\mc O_{\overline \phi}(K\pi)$ and $\overline\phi
\text{-BCC}(K\pi)$ are finite and computable sets. Since $N=\Ker(\pi)$ is computable, then $\overline\phi
\text{-TCC}(K\pi)\pi^{-1}$, $\mc O_{\overline \phi}(K\pi)\pi^{-1}$ and $\overline\phi
\text{-BCC}(K\pi)\pi^{-1}$ are also computable recognizable subsets of $G$, and so  $\overline\phi
\text{-TCC}(K\pi)\pi^{-1}\psi$, $\mc O_{\overline \phi}(K\pi)\pi^{-1}\psi$ and $\overline\phi
\text{-BCC}(K\pi)\pi^{-1}\psi$  are computable  as well.
\qed\\

\begin{corollary}
Let $G$ be a finitely presented group and $\mc C$ be a class of  finite-computable subsets. If $G$ is $\mc C$-generalized twisted conjugacy (resp. Brinkmann's, resp. Brinkmann's conjugacy) separable, then the $\mc C$-genralized twisted conjugacy (resp. Brinkmann's, resp. Brinkmann's conjugacy) problem is decidable.
\end{corollary}

\section{Finite  extensions}
The main goal of this section is to relate twisted conjugacy separability with respect to virtually inner automorphisms in a group $G$ with conjugacy separability of a finite extension of $G$. 

We start by studying the closure of rational subsets of finite extensions of groups, which yields an alternative proof for the fact that finite extensions of product separable groups are product separable \cite{[Cou01]}.

\subsection{Rational subsets}
We know from Theorem \ref{grunschlag-silva-carvalho} that if $H\leq_{f.i.} G$, and $G=Hb_1\cup\cdots\cup Hb_n$, then the rational subsets of $G$ are precisely the subsets of the form $\bigcup_{i\in [n]} L_ib_i$ such that $L_i\in \Rat(H)$. We show the analogous result for separable rational subsets. It is clear that if $H\leq_{f.i} G$, then the inclusion $\iota:H\to G$, when $G$ and $H$ are endowed with the respective profinite topologies, is an open map, as the basis elements for the topology in $H$ are also open in $G$. 
\begin{lemma} \label{separable fi}
Let $G$ be a group, $H\leq_{f.i.} G$ be a finite index subgroup such that $G=Hb_1\cup \cdots \cup Hb_n$ is a decomposition of $G$ as a disjoint union such that $b_1\in H$. A subset $K\subseteq H $ is separable in $G$ if and only if it is separable in $H$.
\end{lemma}
\noindent\textit{Proof.} Since $H$ has finite index, there is a finite index normal subgroup $F\trianglelefteq_{f.i.} H$.

Suppose that $K\subseteq H$ is separable in $G$. Then, for each $g\in G\setminus K$, there is a finite index subgroup $N\trianglelefteq_{f.i.} G$ such that $g\not\in SN$. In particular, for each element $h\in H\setminus K$, there is such an $N$. But then $N\cap F\trianglelefteq_{f.i.} H$ and $S(N\cap F)\subseteq SN$, so  $h\not\in S(N\cap F)$.

Conversely, suppose that $K$ is separable in $H$. Then, $H\setminus K$ is an open set of $H$. But $G\setminus K= H\setminus K \cup Hb_2\cup\cdots\cup Hb_n$ is a finite union of open sets of $G$ since cosets of a finite index subgroup are basis elements and $H\setminus K$ is open in $G$ since the inclusion is an open map.
\qed\\

\begin{proposition}
Let $G$ be a group, $H\leq_{f.i.} G$ be a finite index subgroup such that $G$ admits a decomposition has a disjoint union of the form $G=Hb_1\cup \cdots \cup Hb_n$. Then the  separable rational subsets of $G$ are the ones of the form  $\bigcup_{i\in [n]} L_ib_i$, where  $L_i$ is a separable rational subset of $H$. 
\end{proposition}
\noindent\textit{Proof.} Suppose that $K$ is a separable rational subset of $G$. Then, $K=\bigcup_{i\in [n]} L_ib_i$, where $L_ib_i=K\cap Hb_i$ is closed in the profinite topology, as it is the intersection of two closed subsets.   Since right multiplication by $b_i$, $r_{b_i}$, is also continuous, then $L_i=(L_ib_i)r_{b_i}^{-1}$ is also closed.  By Lemma \ref{separable fi}, $L_i$ is a separable rational subset of $H$.

 Conversely, if $L_i$ is a  separable rational subset of $H$, then it is separable in $G$, by Lemma \ref{separable fi}. Again, since multiplication by $b_i^{-1}$ is continuous, then $L_ib_i$ is separable and $K$ is a finite union of separable subsets, hence separable. Thus $$\overline K=\bigcup_{i\in [n]} \overline{L_ib_i}=\bigcup_{i\in [n]} \overline{L_i}b_i.$$
\qed\\

Notice that the notation for the closure $\overline{L_i}$ is not ambiguous as the closure of $L_i$ in $H$ and in $G$ are the same. Indeed, the closure of $L_i$ in $H$ is closed in $G$ by Lemma \ref{separable fi} and so it contains the closure of $L_i$ in $G$. Moreover, since the closure of $L_i$ in $G$ is contained in $H$ (because it is contained in the closure of $L_i$ in $H$), then it is closed in $H$ and it contains the closure of $L_i$ in $H$, thus being the same. 

\begin{corollary}
Let $G$ be a group, $H\leq_{f.i.} G$ be a finite index subgroup such that $G$ admits a decomposition has a disjoint union of the form $G=Hb_1\cup \cdots \cup Hb_n$. Then, the closure of a rational subset $K=\bigcup_{i\in [n]} L_ib_i\in \Rat(G)$, where $L_i\in \Rat(H)$, is given by $\overline K=\bigcup_{i\in [n]} \overline{L_i}b_i$.
\end{corollary}
\noindent\textit{Proof.} It is clear that, for all $i\in[n]$, $\overline{L_ib_i}\subseteq \overline{L_i}b_i$. Also, $\overline{L_i}\subseteq \overline{L_ib_i}b_i^{-1}$, because  $\overline{L_ib_i}b_i^{-1}$ is a closed subset of $G$ containing $L_i$, and so $\overline{L_i}b_i\subseteq \overline{L_ib_i}$.
\qed\\

A group is said to be \emph{product separable} if, for finitely generated subgroups $H_1,\ldots, H_k$, the product $H_1\cdots H_k$ is closed. Pin an Reutenauer conjectured \cite{[PR91]} that free groups were product separable and proved the following:
\begin{theorem}[Pin-Reutenauer,\cite{[PR91]}]
 If G is a product separable group then the class of all separable rational subsets of G is precisely  the class of subsets of the form $H_1\cdots H_kg$, for some finitely generated subgroups $H_1,\ldots, H_k$ and $g\in G$.
\end{theorem}
 They also showed that product separability of the free group implies Rhodes's type II conjecture (see \cite{[Pin89]}).  Ribes and Zalesskii  \cite{[RZ93]} later proved the conjecture, establishing the result. The class of product separable groups is known to be closed under taking subgroups, finite extensions and free products \cite{[Cou01]}.

From the discussion above, we can provide an alternative proof that the class of product separable groups is
closed under taking finite extensions.
\begin{corollary}
The class of product separable groups is closed under taking finite extensions.
\end{corollary}
\noindent\textit{Proof.} Let $G$ be a group and $H\leq_{f.i.} G$ be a finite index subgroup such that $G$ admits a decomposition has a disjoint union of the form $G=Hb_1\cup \cdots \cup Hb_n$. Suppose that $H$ is product separable and let $K\in \Rat(G)$ be closed. Then, there are rational subsets $L_i\in \Rat(H)$ such that  
$$\bigcup_{i\in [n]} {L_i}b_i=K=\overline K=\bigcup_{i\in [n]} \overline{L_i}b_i.$$
Since the subsets are disjoint, we deduce that $L_i=\overline{L_i}$, and so $L_i$ is closed for all $i\in [n]$. Since $H$ is product separable, then $L_i$ must be of the form $$L_i=\bigcup_{j\in [m_i]} H_{j_1}\cdots H_{j_{s_j}}g_j$$ and so 
$$K=\bigcup_{i\in[n]} L_ib_i=\bigcup_{i\in[n]}\bigcup_{j\in [m_i]} H_{j_1}\cdots H_{j_{s_j}}g_jb_i.$$
Therefore, $G$ is product separable.
\qed\\
\subsection{(Twisted) conjugacy separability}

It is proved in \cite{[Gor86]} that conjugacy separability is not preserved by taking finite extensions and it is not known if twisted conjugacy separability is preserved by taking finite extensions (see \cite{[Ter24]}).  In \cite{[CZ09]}, an example of a non-conjugacy separable group with a finite index conjugacy separable subgroup is presented and a sufficient condition for this property to be preserved by taking finite index subgroups is shown.
In this section, we prove that a finite extension of a group $G$ is (generalized) twisted conjugacy separable group with respect to virtually inner automorphisms if and only if $G$ is (generalized) conjugacy separable.

\begin{theorem} \label{thm: finite extensions}
Let $G$ be a group, $N\trianglelefteq_{f.i.} G$, so that $G$ can be written as a disjoint union $G=Nb_1\cup \cdots \cup Nb_m$, and $K\subseteq G$. For $i\in[m]$, write $L_i=(K\cap Nb_i)b_i^{-1}$.
 If, for all $i\in [m]$, and all $\varphi\in \Via(N)$, the subset $\varphi\text{-}\TCC_N(L_i)$ is separable (in $N$), then  $\alpha(K)$ is separable (in $G$).
 \end{theorem}
\noindent\textit{Proof.} Let $\varphi_i:N\to N$ be such that $x\varphi_i=b_i^{-1}xb_i$. Then $\varphi_i$ is a virtually inner automorphism of $N$: it is injective, as it is the restriction of an automorphism of $G$; it is surjective, since, for $y\in N$, by normality we have $b_iyb_i^{-1}\in N$ and $(b_iyb_i^{-1})\varphi_i=y$; and it is virtually inner, since $\faktor{G}{N}$ is finite, thus $[b_i]$ has finite order, say $p$, in the quotient and $\varphi_i^p$ represents conjugation by $b_i^p$, an element of $N$.

We have that $\alpha(K)=\bigcup_{i=1}^m \alpha(L_ib_i)$. Then, 
$$\alpha(K)=\alpha(K)\cap G=\bigcup_{i=1}^m (\alpha(K)\cap Nb_i)=\bigcup_{i=1}^m\bigcup_{j=1}^m  (\alpha(L_jb_j)\cap Nb_i).$$ We will prove that, for all $i,j\in [m]$, $\alpha(L_jb_j)\cap Nb_i$ is closed in the profinite topology of $G$, which implies that so is $\alpha(K)$.

Put $\beta_{i,j}=\{k\in [m]: b_k^{-1}b_jb_k\in Nb_i\}$ and let $\nu_{k,j}\in N$ be such that $b_k^{-1}b_jb_k=\nu_{k,j}b_r$, for some $r\in [m]$.

We will now prove that $$\alpha(L_jb_j)\cap Nb_i=\bigcup_{k\in \beta_{i,j}} \left(\varphi_j\text{-}\TCC_N(L_j)\right)\varphi_k^{-1}\nu_{k,j}b_i.$$
Let $xb_i\in \alpha(L_jb_j)\cap Nb_i$, where $x\in N$. Then there are some $y\in L_j$, $z\in N$ and $k\in [m]$, such that 
$$xb_i=(b_k^{-1}z^{-1})(yb_j)(zb_k)=(z^{-1}y(z\varphi_j))\varphi_k^{-1}b_k^{-1}b_jb_k.$$ This means that $k\in \beta_{i,j}$ and $b_k^{-1}b_jb_k=\nu_{k,j}b_i$. Hence, 
$xb_i\in  \left(\varphi_j\text{-}\TCC_N(L_j)\right)\varphi_k^{-1}\nu_{k,j}b_i$, for some $k\in \beta_{i,j}.$ Conversely, let $k\in \beta_{i,j}$. Then $b_k^{-1}b_jb_k=\nu_{k,j}b_i$. Now, take an element in $\varphi_j\text{-}\TCC_N(L_j)$, that is an element of the form $z^{-1}y(z\varphi_j)$, for some $z\in N$ and $y\in L_j$. We have that 
$$\left(z^{-1}y(z\varphi_j)\right)\varphi_k^{-1}\nu_{k,j}b_i=\left(z^{-1}y(z\varphi_j)\right)\varphi_k^{-1}b_k^{-1}b_jb_k=(b_k^{-1}z^{-1})(yb_j)(zb_k)\in \alpha(L_jb_j)\cap Nb_i.$$
We know that $\varphi_j\text{-}\TCC_N(L_j)$ is a closed subset of $N$. Since the preimage of a coset of a finite index subgroup by an automorphism is again a coset of a finite index subgroup, then $\varphi_k$ is a continuous automorphism of $N$, when endowed with the profinite topology, so $\left(\varphi_j\text{-}\TCC_N(L_j)\right)\varphi_k^{-1}$ is again closed. Since multiplication by $\nu_{k,j}^{-1}$ is also continuous, and so  $\left(\varphi_j\text{-}\TCC_N(L_j)\right)\varphi_k^{-1}\nu_{k,j}$ is again a closed subset of $N$. By Lemma  \ref{separable fi}, it is closed in $G$ and by continuity of multiplication by an element of $G$, we get that $\left(\varphi_j\text{-}\TCC_N(L_j)\right)\varphi_k^{-1}\nu_{k,j}b_i$ is closed, and then so is 
$$\bigcup_{k\in \beta_{i,j}} \left(\varphi_j\text{-}\TCC_N(L_j)\right)\varphi_k^{-1}\nu_{k,j}b_i=\alpha(L_jb_j)\cap Nb_i.$$
\qed\\

\begin{corollary}
Let $\mc C$ be a class of groups closed under taking finite index subgroups and $\mc K\in \{\Fin,\Coset,\Rat,\Rec,\Alg,\CF\}$ be a class of subsets. If all groups in $\mc C$ are $\mc K$-generalized twisted conjugacy separable with respect to virtually inner automorphisms and $G$ is a group having a finite index subgroup belonging to $\mc C$, then $G$ is $\mc K$-generalized conjugacy separable.
\end{corollary}
\noindent\textit{Proof.}  Let $F\in \mc C$ such that $F\leq_{f.i} G$. Then there is a finite index normal subgroup $N\trianglelefteq_{f.i.}  G$ contained in $F$. Since $\mc C$ is closed under taking finite index subgroups, then $N\in \mc C$.

Put $G=Nb_1\cup\cdots \cup Nb_m$
Let $K\in \mc K(G)$ and put $L_i=(K\cap Nb_i)b_i^{-1}$. If $\mc K=\Fin$, it is obvious that $L_i\in \Fin(N)$ and it follows from the fact that the intersection of two cosets is either empty or a coset of the intersection that if $\mc K=\Coset$, then $L_i\in\Coset(N)$. 
Otherwise, it follows from Theorem \ref{grunschlag-silva-carvalho}, that  $L_i\in \mc K(N)$. Since $N\in\mc C$, then, $\varphi\text{-}\TCC(L_i)$ is separable (in $N$) for all virtually inner automorphisms $\varphi\in \Via(N)$, and it follows from Theorem \ref{thm: finite extensions} that $\alpha(K)$ is separable in $G$. 
\qed\\

Now we will study the converse implication. Let $G=\langle A\mid R\rangle$ be a finitely generated group and $\phi\in \Via(G)$ be such that $\phi^{p_\phi}=\lambda_{\Delta_\phi}$ is the inner automorphism defined by $x\phi^{p_{\phi}}=\dphi^{-1}x\dphi$. We define 
\begin{align*}
G^\phi=\langle A,t \mid R,\, t^{-1}a_it= a_i\phi,\, t^{p_\phi}=\dphi\rangle.
\end{align*}

When the automorphism $\phi$ is clearly set, we will simply write $\Delta$ and $p$ instead of $\dphi$ and $p_\phi$.

\begin{lemma}\label{fi gphi}
Let $G$ be a finitely generated group and $\phi\in \Via(G)$. Then $G\trianglelefteq G^\phi$ and $\left|\faktor{G^\phi}{G}\right|=p_\phi$.
\end{lemma}
\noindent\textit{Proof.} We start by showing that every element of $G^\phi$ can be written in a unique way as $t^kg$ where $g\in G$ and $0<k<p_\phi$. Using the fact that $a_it=t(a_i\phi)$ and $a_it^{-1}=t^{-1}(a_i\phi^{-1})$, given a word in $\widetilde{A\cup\{t\}}$, we can move all the $t's$ to the left and obtain an equivalent word of the form $t^kg$ where $k\in \Z$ and $g\in G$. Since $t^{p_\phi}=\Delta$, we can write $k=r p_\phi +s$ where $r\in \Z$ and $0\leq s<p_\phi$, from where it follows that $t^k=t^s \Delta_\phi^r$
and $t^kg=t^s(\Delta^rg)$. To prove uniqueness, it suffices to observe that if $t^kg=t^\ell h$, then $t^{k-\ell}\in G$, and so $k=\ell$ (notice that $|k-\ell|<p_\phi$)  and $g=h$.

It is clear that $G\trianglelefteq G_\phi$, since, for $g\in G$, we have that $tgt^{-1}=g\phi^{-1}\in G$ and $t^{-1}gt=g\phi\in G$. Also, for elements (written in normal form) $t^kg, t^\ell h\in G^\phi$, we have that  $(t^kg) G=(t^\ell h)G$ if and only if $k=\ell$, hence   $\left|\faktor{G^\phi}{G}\right|=p_\phi$.
\qed\\

\begin{theorem} \label{fi subgroups}
Let $G$ be a finitely generated group, $\phi\in \Via(G)$ and $\mc K\in \{\Fin,\Coset,\Rec,\CF\}$ be a class of subsets. If $G^\phi$ is $\mc K$-generalized conjugacy separable, then  $G$ is $\mc K$-generalized $\phi$-twisted conjugacy separable.
\end{theorem}
\noindent\textit{Proof.} Suppose that  $G^\phi$ is $\mc K$-generalized conjugacy separable and take $K\in \mc K(G)$. We want to prove that  $\phi\text{-}\TCC(K)$ is separable. 

Let $g\in G$.
 For $x\in G$ and $0\leq n< p_\phi$, we have that:
\begin{align*}
(t^nx)(t^{-1}(g\phi^{-1}))(x^{-1}t^{-n})=(t^nx)\cdot g \cdot (x^{-1}\phi)\cdot t^{-n-1}=(xg(x^{-1}\phi))\phi^{-n}t^{-1}.
\end{align*}
Hence, 
\begin{align*}  
\bigcup_{i=0}^{p_\phi-1}(\phitcc_G(g))\phi^{-i}=\alpha_{G^\phi}(t^{-1}(g\phi^{-1}))t.
\end{align*}
Now, we observe that $(\phitcc_G(g))\phi^{-1}\subseteq \phitcc_G(g)$. Indeed, for $x\in G$, 
$$(xg(x^{-1}\phi))\phi^{-1}=(xg)\phi^{-1}x=(xg)\phi^{-1}g(g^{-1}x)\in\phitcc_G(g).$$
Thus, $$\phitcc_G(g)=\bigcup_{i=0}^{p_\phi-1}(\phitcc_G(g))\phi^{-i}=\alpha_{G^\phi}(t^{-1}(g\phi^{-1}))t.$$

Therefore, $$\phitcc_G(K)=\bigcup_{g\in K}\phitcc_G(g)=\bigcup_{g\in K}\alpha_{G^\phi}(t^{-1}(g\phi^{-1}))t=\alpha_{G^\phi}(t^{-1}K\phi^{-1})t.$$
Since $K\in \mc K(G)$ and $G\leq_{f.i.} G^\phi$ by Lemma \ref{fi gphi}, then $K\phi^{-1}\in \mc K(G)$ and  $K\phi^{-1}\in \mc K(G^\phi)$ (see \cite[Lemma 3.4 and Proposition 4.1]{[Car23c]}). By Lemma \ref{herbst letra},  $t^{-1}(K\phi^{-1})\in \mc K(G^\phi)$. Since $G^\phi$ is $\mc K$-generalized conjugacy separable, then $\alpha_{G^\phi}(t^{-1}(K\phi^{-1}))$ is separable in $G^\phi$ and by continuity of multiplication by $t^{-1}$, $\alpha_{G^\phi}(t^{-1}(K\phi^{-1}))t$ is also separable in $G^\phi$. Since it is contained in $G$ and $G\leq_{f.i.} G^\phi$, by Lemma \ref{separable fi}, it is also separable in $G$.
\qed\\

Combining Theorems \ref{thm: finite extensions} and \ref{fi subgroups} in the case where $\mc K$ is the class of singletons, we obtain the following corollary.

\begin{corollary}\label{simple finite equiv}
Let $\mc C$ be a class of finitely generated groups and $\mc{FEC}$ be the class of finite extensions of groups in $\mc C$. Then, groups in $\mc C$ are $\Via$-twisted conjugacy separable if and only if groups in $\mc{FEC}$ are conjugacy separable.
\end{corollary}
Since finitely generated virtually free groups are conjugacy separable \cite{[Dye79],[Ste70]}, we get the following:
\begin{corollary}\label{free via twisted}
Free groups are $\Via$-twisted conjugacy separable.
\end{corollary}

We remark that it is not known if free groups are twisted conjugacy separable and that the previous corollary provides, in particular,  an alternative proof for the main result of \cite{[BBV05]}, where it is proved that the twisted conjugacy problem is decidable for virtually inner automorphisms of the free group. This has been, since then, significantly generalized and the twisted conjugacy problem is known to be decidable for any endomorphism of a finitely generated virtually free group \cite[Theorem 5.4]{[Car23b]}.

Virtually inner automorphisms are a strong restriction of the whole group of automorphisms of a group. For instance, it is easy to see that conjugacy separable groups are $\Via$-Brinkmann conjugacy separable.
\begin{proposition}
Let $G$ be a conjugacy separable group. Then $G$ is  $\Via$-Brinkmann conjugacy separable.
\end{proposition}
\noindent\textit{Proof.} Let $\phi\in \Via(G)$, $x\in G$ and $k>0$ be such that $\phi^k\in \Inn(G)$. Then, if $i\equiv_k j$, we have that $\alpha(x\phi^i)=\alpha(x\phi^j)$, and so  $$\BCC(x)=\bigcup_{i=1}^\infty \alpha(x\phi^i)=\bigcup_{i=1}^k \alpha(x\phi^i)$$ is a finite union of closed sets, thus closed.
\qed\\

We will now prove that free-abelian groups are Coset-generalized twisted conjugacy separable, which yields Coset-generalized conjugacy separability for virtually abelian groups. We  state the result for automorphisms, but the analogous statement holds for arbitrary endomorphisms.
\begin{proposition}
	Free-abelian groups are Coset-generalized twisted conjugacy separable, that is, if $G$ is a free-abelian group,  $\phi\in\Aut(G)$ is an automorphism,  and $Hx$ is  a coset of a finitely generated subgroup $H\leq_{f.g.} G$, the set of all $\phi$-twisted conjugates of elements in $Hx$, $\phitcc(Hx)$, is  separable.\end{proposition}
\noindent\textit{Proof.} Let $H\leq_{f.g.} \Z^m$, $x\in \Z^m$ and $M$ be the matrix representing the automorphism $\phi$. Then we have:
\begin{align*}
\phitcc(Hx)&=\bigcup_{z\in\Z^m} (z^{-1}\phi)Hx z\\
&=\bigcup_{z\in\Z^m} Hx (z(\Img(I-M)))\\
&=H\Img(I-M)x.
\end{align*}
Thus, $\phitcc(Hx)$ is itself a coset of a subgroup of $\Z^m$. Since cosets of subgroups of free-abelian groups are separable \cite[Corollary 4.3]{[Del98]}, then so is $\phitcc(Hx)$.
\qed\\

Virtually abelian groups include metacyclic groups and locally finite linear groups over characteristic zero.
\begin{corollary}\label{virtually abelian conjugacy}
Every virtually abelian group is $\Coset$-generalized conjugacy separable.
\end{corollary}

\begin{remark}
Notice that since not all rational subsets of free-abelian groups are separable and conjugacy in abelian groups coincides with equality, then 
free-abelian groups are not $\Rat$-generalized conjugacy separable. However, in free-abelian groups, given an automorphism $\phi$ and a rational subset $K$, it is decidable whether $\phitcc(K)$ is closed, since $\phitcc(K)$ is rational \cite{[Del98]}.
\end{remark}

\subsection{Virtually free times free groups}

We will prove that virtually free times free groups are conjugacy separable by proving that free times free groups are $\Via$-twisted conjugacy separable. If both groups have rank 1, then this is already known, as such groups are virtually abelian. We start by dealing with groups of the form $F_n\times \Z$ with $n\geq 2$ and then with groups of the form $F_n\times F_m$ with $n,m \geq 2$.

Given a virtually inner automorphism $\phi\in \Via(G)$, we put $$\mathfrak o (\phi)=\min\{k>0\mid \phi^k\in \Inn(G)\}.$$
\begin{lemma}\label{inner order}
Let $G$ be a group and $\phi\in \Via(G)$. Then $\phi^k\in \Inn(G)$ if and only if $k\in \mathfrak o(\phi)\Z$.
\end{lemma}
\noindent\textit{Proof.} Suppose that $\phi^{\mathfrak o(\phi)}=\lambda_z$, where $\lambda_z$ denotes conjugation by $z$. Then, $\phi^{m\mathfrak o(\phi)}=\lambda_{z^m}\in \Inn(G)$, for all $m\in \Z$. Conversely, suppose that  $\phi^k\in \Inn(G)$ and $k\not\in \mathfrak o(\phi)\Z$. We have that $|k|>\mathfrak o(\phi)$, since $\phi^{|k|}\in \Inn (G)$ and $\mathfrak o(\phi)$ is minimal.
Let $s=\max\{m\in \N\mid m\mathfrak o(\phi) < |k|\}$. Then, $\phi^{|k|-s\mathfrak o(\phi)}\in \Inn(G)$ and $0<|k|-s\mathfrak o(\phi)<\mathfrak o(\phi)$, which contradicts the minimality of $\mathfrak o(\phi)$.
\qed\\

A group $F_n\times \Z^m$ 
admits a presentation of the form: 
\begin{equation*}
\langle x_1,\ldots,x_n,t_1,\ldots,t_m \mid t_it_j \,=\, t_jt_i,\, t_ix_k=x_kt_i \,, i,j\in[1,m], \ k\in[1,n] \rangle,
 \end{equation*}
where $n,m \geq 0$.
Note that, given a word in the generators, we can use the commutativity relations to move all the $t_j$'s  to the left and describe the elements in $G$ by a normal form $ t_1^{a_1} \cdots t_m^{a_m}\, u(x_1,\ldots,x_n)$, which we abbreviate as $t^au$, where $u \in F_n$ and ${a} = (a_1,\ldots,a_m) \in \Z^m$.

In \cite{[DV13]}, it is shown that automorphisms of free-abelian times free groups $F_n\times \Z^m$ are of the form $$t^au\mapsto t^{aQ+\au P}u\phi,$$ where $\phi\in \Aut(F_n)$, $Q\in \GL_m(\Z)$, $P\in \mc M_{n\times m}(\Z)$ and $\au$ denotes the abelianization of $u$. Such an automorphism is completely characterized by the matrices $Q,\, P$ and the automorphism $\phi$, thus being denoted by $\Psi_{\phi,Q, P}$.

Also, \cite[Proposition 1.4]{[CD24c]} shows that 
\begin{align}
\label{powers aut}
\Psi_{\phi, Q, P}^k=\Psi_{\phi^k,Q^k,P^{(k)}},
\end{align} where $\phi^{ab}:\Z^n\to \Z^n$ is the abelianization map of $\phi$ and $P^{(k)}=\sum_{i=1}^k (\phi^{ab})^{i-1}PQ^{k-i}$.

We will now describe virtually inner automorphisms of free-abelian times free groups, thus describing, in particular, virtually inner automorphisms of $F_n\times \Z$.

\begin{lemma} \label{virtinner fatf}
Virtually inner automorphisms of a free-abelian times free group $F_n\times \Z^m$ are precisely the automorphisms $\Psi_{\phi, Q, P}$ where $\phi\in \Via(F)$, $Q$ has finite order and there is some $k\in  \mathfrak o(\phi)\N$ such that 
$$\sum_{i=1}^k(\au\phi^{ab})^{i-1}PQ^{k-i}=0$$
for all $u\in F_n$.
\end{lemma}
\noindent\textit{Proof.} Suppose that $\Psi_{\phi, Q, P}$ is such that $\Psi^k=\lambda_{t^by}$, where $k\in \N$ and $\lambda_{t^b y}$ denotes the inner automorphism defined by conjugation by $t^by\in F_n\times \Z^m$. Then, by (\ref{powers aut}),
$$t^a(y^{-1}uy)=(t^au)\lambda_{t^by}=(t^{a}u)\Psi_{\phi, Q, P}^k=(t^{a}u)\Psi_{\phi^k,Q^k,P^{(k)}}=t^{aQ^k+\au P^{(k)}}(u\phi^k)$$
for all $t^au\in F_n\times \Z^m$. Hence, for all $u\in F_n$, $u\phi^k=y^{-1}uy$, that is $\phi\in \Via(F_n)$ as $\phi^k$ represents conjugation by $y$ and, by Lemma \ref{inner order}, $k\in \mathfrak o(\phi)\N$; putting $u=1$, we obtain that for all $a\in \Z^m$, $aQ^k=a$, so $Q^k=I_m$ and $Q$ has finite order; and for all $u\in F_n$ and $a\in \Z^m$,  $aQ^k+\au P^{(k)}=a$, which implies that $\au P^{(k)}=0$, since $aQ^k=a$.

Now, let $\Psi_{\phi, Q, P}$ be such that $\phi^{k_1}$ is inner (and represents conjugation by $y$) for some $k_1\in \N$, $Q^{k_2}={I_m}$ for some $k_2\in \N$,  and there is some $s\in \N$ such that 
$\sum_{i=1}^{sk_1}(\au\phi^{ab})^{i-1}PQ^{sk_1-i}=0$
for all $u\in F_n$.
Then, for all $u\in F_n$
$$u\Psi_{\phi, Q, P}^{sk_1}=u\Psi_{\phi^{sk_1}, Q^{sk_1}, P^{(sk_1)}}=t^{\au P^{(sk_1)}}u\phi^{sk_1}=u\phi^{sk_1}=y^{-s}uy^s,$$ thus
$$u\Psi_{\phi, Q, P}^{sk_1k_2}=y^{-sk_2}uy^{sk_2}$$
and, for all $a\in \Z^m$,
$$t^a\Psi_{\phi, Q, P}^{sk_1k_2}=t^{aQ^{sk_1k_2}}=t^a.$$
Hence, $$(t^au)\Psi_{\phi, Q, P}^{sk_1k_2}=t^ay^{-sk_2}uy^{sk_2},$$
which means that $\Psi_{\phi, Q, P}^{sk_1k_2}$ is inner and represents conjugation by $y^{sk_2}$.
\qed\\

We will now prove two auxiliary lemmas that will be useful to proving the main result of this section.

\begin{lemma}\label{1chega}
Let $G$ be a group. Then $G$ is $\Via$-twisted conjugacy separable if and only if $\phitcc(1)$ is separable for all $\phi\in \Via(G)$.
\end{lemma}
\noindent\textit{Proof.}
Since, for all $x,g\in G$, 
\begin{align}\label{tcc 1 chega}
(x^{-1}\phi) gx=gg^{-1}(x^{-1}\phi) gx=g\left((x^{-1}\phi\lambda g) x\right),
\end{align} it follows that, for all $g\in G$,  
$\phitcc(g)=g(\phi\lambda_g\text{-TCC($1$)})$. Since $\phi\lambda_g\in \Via(G)$ by \cite[Proposition 2.5]{[LS11]}, $(\phi\lambda_g\text{-TCC($1$)})$ is separable, which, by continuity of left multiplication, implies that so is $\phitcc(g)$.
\qed\\

Given an element $x\in G$, $\phi\in \Aut(G)$ and a subset $K\subseteq G$, we define $\phitcc(x,K)=\{(y^{-1}\phi)xy\mid y\in K\}$ as the set of twisted conjugates of $x$ with a twisted conjugator in $K$. Similar notions could be defined for the remaining problems (conjugacy, Brinkmann and Brinkmann conjugacy). However, only the twisted conjugacy case will be relevant to us. We remark that these are the separability notions corresponding to the decidability problems with constraints on the set of solutions (see, for example, \cite{[LS11]}).

\begin{lemma} \label{free fi constraints}
Let $F$ be a free group, $\phi\in \Via(F)$ and $H\leq_{f.i.} F$. Then, the set $\phitcc(1,Hz)=\{(x^{-1}\phi)x\mid x\in Hz\}$ is separable for all $z\in F$.
\end{lemma}
\noindent\textit{Proof.} By Lemma \ref{fully invariant}, there is a fully invariant finite index subgroup $H'\leq_{f.i.} H$. Since $H'$ is fully invariant, then $\phi|_{H'}\in \Aut(H')$. Also, since $H'$ is fully invariant, $H'\trianglelefteq F$, as it is invariant under the action of inner automorphisms. Moreover, $\phi|_{H'}\in \Via(H')$: since $\phi\in \Via(F)$, then there is some $k\in \N$ such that $\phi^k=\lambda_z$ represents conjugation by $z\in F$ and there must be some $s>0$ such that $z^s\in H'$ since $z$ has finite order in the quotient $\faktor{F}{H'}$, hence $(\phi|_{H'})^{sk}=\lambda_{z^{sk}}\in \Inn(H')$.
We have that $H=H'b_1\cup\cdots\cup H' b_m$, so  
\begin{align*}
\phitcc(1,Hz)&=\bigcup_{i=1}^m \{(x^{-1}\phi)x\mid x\in H'b_iz\}\\
&=\bigcup_{i=1}^m \{(z^{-1}b_i^{-1}x^{-1})\phi xb_iz\mid x\in H'\}\\
&=\bigcup_{i=1}^m \left((z^{-1}b_i^{-1})\phi\right)\{(x^{-1}\phi)x\mid x\in H'\}\,b_iz\\
&=\bigcup_{i=1}^m \left((z^{-1}b_i^{-1})\phi\right) \left(\phi|_{H'}\text{-TCC$_{H'}$}(1)\right)\,b_iz
\end{align*}
By Corollary \ref{free via twisted}, $\phi|_{H'}\text{-TCC$_{H'}$}(1)$ is separable in $H'$, and by Lemma \ref{separable fi}, it follows that it is separable in $F$. By continuity of left and right multiplication,  $(z^{-1}b_i^{-1})\phi \left(\phi|_{H'}\text{-TCC$_{H'}$}(1)\right)b_iz$ is separable for all $i\in [m]$ and so is $\phitcc(1,Hz)$.
\qed\\

We can now prove that virtually free times infinite cyclic groups are conjugacy separable. We remark that these include non-elementary unimodular generalized Baumslag-Solitar groups (see \cite[Proposition 2.6]{[Lev07]}) and that, among free-by-cyclic group with a polynomially growing monodromy, these are precisely those which are subgroup separable \cite{[Kud24]}.

\begin{theorem}
Virtually $F_n\times \Z$ groups are conjugacy separable.
\end{theorem}
\noindent\textit{Proof.} In view of Theorem \ref{fi subgroups}, it suffices to prove that $F_n\times \Z$ is $\Via$-twisted conjugacy separable.
Let $\Psi\in \Via(F_n\times \Z)$. Since $\Z$ has only two automorphisms, from Lemma \ref{virtinner fatf}, it follows that there are some $\phi\in \Via(F)$ and $P\in \Z^n$ for which there is some $k\in  \mathfrak o(\phi)\N$ such that 
 $\Psi$ is  defined as
\begin{align} \label{auto1}
t^mu\mapsto t^{m+\au P}u\phi,
\end{align}
with
$$\sum_{i=1}^k(\au\phi^{ab})^{i-1}P=0,$$ for all $u\in F_n$,
or $\Psi$ is  defined as
\begin{align}\label{auto2}
t^mu\mapsto t^{-m+\au P}u\phi,
\end{align}
with 
$$\sum_{i=1}^k(-1)^{k-i}(\au\phi^{ab})^{i-1}P=0,$$
for all $u\in F_n$.

 By Lemma \ref{1chega}, we only have to prove that $\psitcc(1)$ is closed.

Assume that $\Psi$ is as in (\ref{auto1}). Notice that
$$((u^{-1}t^{-p})\Psi) t^pu=t^{-\au P}(u^{-1}\phi)u.$$
Let $t^qv\not\in\psitcc(1)$. If $v\not\in\phitcc(1),$ then, by Corollary \ref{free via twisted}, there is a homomorphism  $\xi:F_n\to G$ onto a finite group $G$ such that  $v\xi\not\in \phitcc(1)\xi$. Thus, putting $\xi':F_n\times \Z\to G$ as being defined by $t^mx\mapsto x\xi$, we get that $(t^q v)\xi'=v\xi\not\in \psitcc(1)\xi'$, since $\psitcc(1)\xi'=\phitcc(1)\xi$.

If, on the other hand, $v\in \phitcc(1)$, then $t^qv=t^q(y^{-1}\phi)y$ is such that $q\neq -\overrightarrow{y} P$.
It is easy to see that the mapping $\zeta:F_n\times \Z\to \Z$ defined by $$t^a u\mapsto -ka-\sum\limits_{i=1}^{k-1}(k-i)\overrightarrow{u\phi^{i-1}}P$$ is a homomorphism. Also, for $t^{-\ax P}(x^{-1}\phi) x\in \psitcc(1),$ we have that 
\begin{align*}
(t^{-\ax P}x^{-1}(\phi)x)\zeta&= k\ax P-\sum\limits_{i=1}^{k-1}(k-i)\overrightarrow{(x^{-1}\phi^i)(x\phi^{i-1})}P\\
&= k\ax P-\sum\limits_{i=1}^{k-1}(k-i)\overrightarrow{(x^{-1}\phi^i)}P -\sum\limits_{i=1}^{k-1}(k-i)\overrightarrow{(x\phi^{i-1})}P\\
&= k\ax P+\sum\limits_{i=1}^{k-1}(k-i)\overrightarrow{(x\phi^i)}P -\sum\limits_{i=1}^{k-1}(k-i)\overrightarrow{(x\phi^{i-1})}P\\
&= k\ax P-(k-1)\ax P+\sum\limits_{i=1}^{k-1}\overrightarrow{(x\phi^i)}P \\
&=0,
\end{align*}
and so $\psitcc(1)\zeta=\{0\}$.
Similarly, we see that 
$$(t^q(y^{-1}\phi)y)\zeta = -kq  -(k-1)\ay P+\sum\limits_{i=1}^{k-1}\overrightarrow{(y\phi^i)}P=-k(q+\ay P)\neq 0.$$
Since $\Z$ is residually finite, there is a homomorphism $\xi:\Z\to G$ onto a finite group $G$ such that $k(q+\ay P)\xi\neq 0\xi$, so $(t^qv)\zeta\xi\not\in \psitcc(1)\zeta\xi$.

Now, we assume that $\Psi$ is as in (\ref{auto2}). Notice that, in this case, 
$$((u^{-1}t^{-p})\Psi) t^pu=t^{2p -\au P}(u^{-1}\phi)u.$$
So, $$t^a x\in \psitcc(1)\iff   \exists u\in F_n \, \begin{cases} x=(u^{-1}\phi)u\\
a\in -\au P+ 2\Z
\end{cases}.$$

Let $t^qv\not\in\psitcc(1)$. If $v\not\in\phitcc(1),$ we can proceed as above using Corollary \ref{free via twisted}. Hence, we may assume that $t^qv=t^q(y^{-1}\phi)y$ and $q\not\in -\ay P +2\Z$.

Let $H=\{u\in F_n\mid \au P \equiv 0\mod 2\}$. Then  $H\trianglelefteq F_n$ and $[F_n:H]=2$, so $F_n=H\cup Hw$ where $\overrightarrow{w}P\equiv 1 \mod 2$. 

By Lemma \ref{free fi constraints}, there is a homomorphism $\xi:F_n \to G$ onto a finite group $G$ such that $(y^{-1}\phi y)\xi\not\in \phitcc(1,H)\xi$ if $y\not \in H$ and there is a homomorphism $\xi:F_n \to G$ onto a finite group $G$ such that $(y^{-1}\phi y)\xi\not\in \phitcc(1,Hw)\xi$ if $y\not \in Hw$.

Consider the homomorphism $\zeta:F_n\times \Z\to \Z_2\times G$ defined by $t^au\mapsto (a\mod 2,  u\xi)$. 
We will prove that $(t^q(y^{-1}\phi)y)\zeta\not\in \psitcc(1)\zeta.$ Consider some element $t^m(z^{-1}\phi) z$ such that $m\in -\overrightarrow{z}P+2\Z$, so that 
 $t^m(z^{-1}\phi) z\in\psitcc(1)$ and suppose that $(t^q(y^{-1}\phi)y)\zeta=(t^m(z^{-1}\phi) z)\zeta$. Then 
 $$\begin{cases}
 m\equiv q \mod 2\\
 ((y^{-1}\phi)y)\xi=((z^{-1}\phi) z)\xi
 \end{cases}.$$
 We can deduce, by the definition of $\xi$ that $Hz= Hy$, that is $\overrightarrow{z}P\equiv \ay P \mod 2$. But then, $q\equiv m\equiv \overrightarrow{z}P\equiv \overrightarrow{y}P \mod 2$, which is absurd.
\qed\\

In \cite{[Car23d]}, automorphisms of $F_n\times F_m$ are described. If $n\neq m$, the only automorphisms of $F_n\times F_m$ are those of the form 
\begin{align}\label{type vi}
(x,y)\mapsto (x\phi,x\psi),
\end{align}
 where $\phi\in \Aut(F_n)$ and $\psi\in \Aut(F_m)$; and if $n=m$, there are also automorphisms of the form 
 \begin{align} \label{type vii}
 (x,y)\mapsto (y\psi,x\phi),
 \end{align} where $\phi,\psi\in \Aut(F_n)$. 
 Following the terminology in \cite{[Car23d]}, we will say that automorphisms of the form (\ref{type vi}) are of type VI and automorphisms of the form (\ref{type vii}) are of type VII.

\begin{theorem}\label{main direct}
Virtually free times free groups are conjugacy separable.
\end{theorem}
\noindent\textit{Proof.} It is clear that automorphisms of type VII cannot be inner: if $\Phi\in \Aut(F_n\times F_n)$ is defined by $(x,y)\mapsto (y\psi,x\phi)$, then, for $x\in F_n\setminus\{1\}$,  $(x,1)\Phi=(1,x\phi)\not\sim (x,1)$. 
Also, powers of automorphisms of type VI are still automorphisms of type VI, but the behavior of type VII automorphisms is different: even powers of type VII automorphisms are automorphisms of type VI while odd powers of automorphisms of type VII are still of type VII.

Let $\Psi\in \Via (F_n \times F_m)$ be a virtually inner automorphism of type VI. Then $\Psi$ is defined by $(x,y)\mapsto (x\phi,y\psi)$ where $\phi\in \Via(F_n)$ and $\psi\in \Via(F_m)$ since $\Inn(F_n\times F_m)\simeq \Inn(F_n)\times \Inn(F_m)$: if $\Psi^k=\lambda_{(x,y)}$, then $\phi^k=\lambda_x$ and $\psi^k=\lambda_y$.

So, in this case $\psitcc((1,1))=\phitcc(1)\times \Psitcc(1)$ is closed as it is a direct product of closed sets in each factor.

Now, suppose that $m=n$ and $\Psi$ is a virtually inner automorphism of type VII of the form $(x,y)\mapsto (y\psi,x\phi)$. In this case, $\mathfrak o(\Psi)$ must be even and  $\psi\phi, \phi\psi\in \Via(F_n)$. Let $\xi\in \Aut(F_n\times F_n)$ be the type VI automorphism defined by $(x,y)\mapsto (x,y\phi^{-1})$ and $m:F_n\times F_n\to F_n$ be the multiplication map (which is continuous), that is, $(x,y)m=xy.$ We will prove that 
\begin{align}\label{twisted type vii}
\psitcc((1,1))\xi=\phi\psi\text{-TCC(1)}m^{-1}.
\end{align}
Let $\left((x^{-1},y^{-1})\Psi\right) (x,y)\in \psitcc((1,1))$. Then, 
\begin{align*}
\left[\left((x^{-1},y^{-1})\Psi\right) (x,y)\right]\xi m &=\left[\left((y^{-1}\psi,x^{-1}\phi)\right) (x,y)\right]\xi m\\
&=\left[\left((y^{-1}\psi)x,(x^{-1}\phi)y\right)\right]\xi m\\
&=\left((y^{-1}\psi)x,x^{-1}(y\phi^{-1})\right) m\\
&=(y^{-1}\psi)xx^{-1}(y\phi^{-1})\\
&=(y^{-1}\psi)(y\phi^{-1})\\
&=\left((y^{-1}\phi^{-1})\phi\psi\right)(y\phi^{-1})\in \phi\psi\text{-TCC(1)}
\end{align*}
Conversely, given $(x,y)\in F_n\times F_n$ such that $xy\in  \phi\psi\text{-TCC(1)},$ there must be some $z\in F_n$ such that  $xy=(z^{-1}\phi\psi)z$, so $x=(z^{-1}\phi\psi)zy^{-1}$  and $(x,y)=((z^{-1}\phi\psi)zy^{-1},y)$, for some $z\in F_n$. Hence 
$$(x,y)\xi^{-1}=((z^{-1}\phi\psi)zy^{-1},y\phi)=\left((yz^{-1},z^{-1}\phi)\right)\Psi (zy^{-1},z\phi)\in \psitcc((1,1)).$$ Since $\phi\psi\text{-TCC(1)}$ is closed in $F_n$ and $m$ is continuous, then $\psitcc((1,1))\xi=\phi\psi\text{-TCC(1)}m^{-1}$ is closed in $F_n\times F_n$.

Therefore, given $(u,v)\not\in \psitcc((1,1))$, we have that $(u,v)\xi\not\in  \phi\psi\text{-TCC(1)}m^{-1}$ and there is a homomorphism $\zeta: F_n\times F_n\to G$ onto a finite group $G$ such that $(u,v)\xi\zeta\not\in\phi\psi\text{-TCC(1)}m^{-1}\zeta =\psitcc((1,1))\xi\zeta$, and $\xi\zeta$ separates $(u,v)$ and $\psitcc(1,1))$.
\qed\\

The twisted conjugacy problem is known to be decidable for free-abelian times free groups \cite[Theorem 3.6]{[CD24c]}.
Looking into the proof of the result above, we can see that the twisted conjugacy problem is decidable for free times free groups: if we are given an automorphism of type VI, the problem reduces to two instances of the twisted conjugacy problem for free groups; if we receive as input two elements $(x,y),(z,w)\in F_n\times F_n$ and a type VII automorphism  $\Psi$ defined by $(u,v)\mapsto (v\psi,u\phi)$, deciding if $(x,y)$ and $(z,w)$ are $\Psi$-twisted conjugate amounts to deciding whether $(z^{-1}x,w^{-1}y)$ and $(1,1)$ are $\Psi\lambda_{(z,w)}$-twisted conjugate,
by (\ref{tcc 1 chega}). Notice that $\Psi\lambda_{(z,w)}$ is the type VII automorphism defined by $(x,y)\mapsto (y\psi\lambda_z,x\phi\lambda_w)$. By (\ref{twisted type vii}), this condition is equivalent to checking if 
$$z^{-1}x(w^{-1}y)\phi^{-1} =(z^{-1}x,w^{-1}y)\xi m\in (\phi\lambda_w\psi\lambda_z)\text{-TCC(1)},$$
which is decidable since the twisted conjugacy problem is decidable for free groups \cite[Theorem 1.5]{[BMMV06]}. Hence, the following is a natural consequence of the proof of the previous theorem.

\begin{corollary}
The twisted conjugacy problem is decidable for the direct product of two free groups. 
\end{corollary}

\section{Cyclic extensions} 
We now turn our focus to infinite cyclic extensions.
Let $A=\{a_1,\ldots, a_n\}$, $G=\langle A\mid R\rangle$ be a group and $\phi\in \Aut(G)$. A $G$-by-$\Z$ group has the form 
\begin{align}
\label{presentation semidirect}
G\rtimes_\phi \Z=\langle A,t \mid R, t^{-1}a_it= a_i\phi\rangle.
\end{align}

Every element of $G\rtimes_\phi \Z$ can be written in a unique way as an element of the form $t^ag$, where $a\in \Z$ and $g\in G$. 

In \cite[Theorem 3.8]{[FT06]}, a connection between conjugacy separability of $G\rtimes \Z$ and twisted conjugacy separability of $G$ is established. Here, we prove not only a connection between conjugacy separability of $G\rtimes \Z$  with Brinkmann conjugacy separability, but also show the result for the generalized versions of the concepts.

\begin{lemma}\label{lem: equivs}
Let $G$ be a group, $x\in G$ and $\phi\in \Aut(G)$.
Then we have that:
\begin{enumerate}
\item $\alpha_{G\rtimes_\phi \Z}(x)=\phibcc_G(x)$;
\item $\alpha_{G\rtimes_\phi \Z}(tx)=t\phitcc(x)$.
\end{enumerate}
\end{lemma}
\noindent\textit{Proof.} Clearly, in $G\rtimes_\phi\Z$, if $t^ax\sim t^by$, then $a=b$.
We start by proving 1:
\begin{align*}
x\sim y &\iff \exists c\in \Z, z\in G : (z^{-1}t^{-c})x(t^cz)=y\\
&\iff  \exists c\in \Z, z\in G: z^{-1}(x\phi^c)z=y\\
&\iff  y\in \phibcc_G(x)
\end{align*}

Now, to prove 2, we have that: 
\begin{align*}
tx\sim ty &\iff \exists c\in \Z, z\in G : (z^{-1}t^{-c})tx(t^cz)=ty\\
&\iff  \exists c\in \Z, z\in G: z^{-1}t(x\phi^c)z=ty\\
&\iff  \exists c\in \Z, z\in G: t(z^{-1}\phi)(x\phi^c)z=ty\\
&\iff  \exists c\in \Z, z\in G: \left[\left(z^{-1}(x\phi^{c-1})\cdots(x\phi) x \right)\phi \right]x  \left(x^{-1}(x^{-1}\phi)\cdots (x^{-1}\phi^{c-1})z\right)=y.\\
&\iff x\sim_\phi y.\\
\end{align*}
\qed\\

\begin{theorem}\label{main generalized descend}
Let $\mc K\in \{\Fin, \Coset,\Rat,\Alg\}$. If the class of $G\rtimes \Z$ groups is $\mc K$-generalized conjugacy separable, then $G$ is $\mc K$-generalized twisted conjugacy separable and $\mc K$-generalized Brinkmann conjugacy separable.
\end{theorem}
\noindent\textit{Proof.} We start by proving that $G$ is $\mc K$-generalized Brinkmann conjugacy separable. Let   $\phi\in \Aut(G)$, $K\in \mc K(G)$. We want to prove that $S=\bigcup_{x\in K} \phibcc_G(x)$ is separable, that is, for all $y\in G\setminus S$, there is a finite group $F$ and a homomorphism $p: G\to F$ such that $p(y)\not\in p(S)$.
Now,  in $G\rtimes_\phi \Z$, we have that $\alpha(K)=S$, by Lemma \ref{lem: equivs}. Since we are assuming that $G\rtimes_\phi \Z$ is $\mc K$-generalized conjugacy separable and $y\not\in \alpha(K)$, there is some finite group $F'$ and a homomorphism $q:G\rtimes_\phi \Z\to F'$ such that $q(y)\not\in q(\alpha(K))=q(S).$ Hence, taking $p=q|_G$, we have that $G$ is $\mc K$-generalized Brinkmann conjugacy separable. 

Now, we deal with the generalized twisted conjugacy separability case. Let   $\phi\in \Aut(G)$, $K\in \mc K(G)$ and $S=\bigcup_{x\in K} \phitcc_G(x)$. We want to prove that $S$ is separable. It is easy to see that $K\in \mc K(G\rtimes_\phi\Z)$, and, by Lemma \ref{herbst letra},  $tK\in \mc K(G\rtimes_\phi\Z)$. Let $T=\alpha_{G\rtimes_\phi \Z}(tK)$  and $y\in G\setminus S$. Then $T=tS$ and so  $ty\not\in T$. 
Hence, there is a finite group $F'$ and a homomorphism $q:G\rtimes_\phi \Z\to F'$ such that $q(ty)\not\in q(T)$. Take $p=q|_G$. If $p(y)\in p(S)$, then this would mean that $q(ty)=q(t)q(y)\in q(t)p(S)=q(tS)=q(T)$, which is absurd. Thus $p(y)\not\in p(S)$.
\qed\\

\begin{corollary}\label{simple cyclic equiv}
If the class of $G\rtimes \Z$ groups is conjugacy separable, then $G$ is twisted conjugacy separable and Brinkmann conjugacy separable.
\end{corollary}

As a corollary, we obtain separability results for virtually polycyclic groups (twisted conjugacy separability was already proved in \cite{[FT06]}).

\begin{corollary}\label{virt polycyclic twisted Brink}
Virtually polycyclic groups are twisted conjugacy separable and Brinkmann conjugacy separable.
\end{corollary}
\noindent\textit{Proof.} It is proved in \cite[Proposition 4.3]{[Car23c]} that, if $G$ is virtually polycyclic and $\phi\in\Aut(G)$, then $G\rtimes_\phi \Z$ is virtually polycyclic. Since polycyclic groups are conjugacy separable, the result follows from Theorem \ref{main generalized descend}.
\qed\\

Since in the case of abelian groups, conjugacy is simply equality, we obtain that the orbit of any element by an automorphism of an abelian group is closed in the profinite topology.
\begin{corollary}
Finitely generated abelian groups are Brinkmann equality separable.
\end{corollary}

When decidability is concerned, for the simple case, we have equivalence, that is, $G$-by $\Z$ groups have decidable conjugacy problem if and only if $G$ has decidable twisted conjugacy and Brinkmann conjugacy problems. We do not know if the same holds for separability.

\begin{question}
Does the converse of Theorem \ref{main generalized descend} hold? Does it hold for particular classes of subsets (for example, singleton sets)?
\end{question}

It is not known whether free-by-cyclic groups are conjugacy separable or not. It is then a natural question to ask if free groups are Brinkmann (conjugacy) separable.

\begin{question}
Are free groups of finite rank Brinkmann (conjugacy) separable?
\end{question}

\section*{Acknowledgements}
The author was partially supported by
CMUP, member of LASI, which is financed by national funds through FCT – Fundação
para a Ciência e a Tecnologia, I.P., under the projects with reference UIDB/00144/2020
and UIDP/00144/2020.

\bibliographystyle{plain}
\bibliography{Bibliografia}

 \end{document}